\documentclass[12pt,twoside,a4paper]{amsart}
\usepackage{amsmath,amssymb,amscd,amsfonts,latexsym,a4wide}

\numberwithin{equation}{section}

\newcommand{\Z}{{\mathbb Z}}
\newcommand{\Q}{{\mathbb Q}}

\newcommand{\Ga}{\Gamma}

\newcommand{\bx}{{\bf x}}
\newcommand{\by}{{\bf y}}

\newcommand{\bu}{{\bf u}}
\newcommand{\ba}{{\bf a}}
\newcommand{\bb}{{\bf b}}

\newcommand{\bone}{{\bf 1}}

\newcommand{\kdots}{,\ldots,}
\newcommand{\Kn}{(K^*)^n}

\newcommand{\rank}{\mbox{$\mathrm{rank}\,$}}

\renewcommand{\gcd}{\mbox{$\mathrm{gcd}$}}

\newcommand{\eps}{\varepsilon}

\newcommand{\s}{\sigma}
\newcommand{\bus}{{\bf u}_{\sigma}}
\newcommand{\be}{\begin{equation}}
\newcommand{\ee}{\end{equation}}
\newcommand{\formref}[1]{(\ref{#1})}
\newcommand{\romref}[1]{{\rm (\ref{#1})}}

\renewcommand{\leq}{\leqslant}
\renewcommand{\geq}{\geqslant}

\title[Linear equations with unknowns from a multiplicative group]{Linear equations
with unknowns from a multiplicative group whose solutions lie in a
small number of subspaces}
\author[J.-H.~EVERTSE]{Jan-Hendrik~EVERTSE}
\address{Universiteit Leiden, Mathematisch Instituut, Postbus 9512,
NL-2300 RA Leiden}
\email{evertse@math.leidenuniv.nl}

\begin{document}

\begin{abstract}
Let $K$ be a field of characteristic $0$ and let $\Kn$ denote the
$n$-fold cartesian product of $K^*$, endowed with coordinatewise
multiplication. Let $\Ga$ be a subgroup of $\Kn$ of finite rank.
We consider equations (*) $a_1x_1+\cdots+a_nx_n=1$ in ${\bf
x}=(x_1\kdots x_n)\in\Ga$, where ${\bf a}=(a_1\kdots a_n)\in \Kn$.
Two tuples $\ba ,\bb\in\Kn$ are called $\Ga$-equivalent if there
is a $\bu\in \Ga$ such that $\bb =\bu\cdot \ba$. Gy\H{o}ry and the
author \cite{EG} showed that for all but finitely many
$\Ga$-equivalence classes of tuples $\ba\in\Kn$, the set of
solutions of (*) is contained in the union of not more than
$2^{(n+1)!}$ proper linear subspaces of $K^n$. Later, this was
improved by the author \cite{Evb} to $(n!)^{2n+2}$. In the present
paper we will show that for all but finitely many
$\Ga$-equivalence classes of tuples of coefficients, the set of
non-degenerate solutions of (*) (i.e., with non-vanishing subsums)
is contained in the union of not
more than $2^n$ proper linear subspaces of $K^n$. Further we give
an example showing that $2^n$ cannot be replaced by a quantity
smaller than $n$.
\\[0.3cm]
\\[0.2cm]
\\[0.3cm]
\end{abstract}

\maketitle

\section{Introduction}

Let $K$ be a field of characteristic $0$. Denote by $(K^*)^n$ the
$n$-fold direct product of the multiplicative group $K^*$. The
group operation of $\Kn$ is coordinatewise multiplication, i.e.,
if $\bx =(x_1\kdots x_n)$, $\by =(y_1\kdots y_n)\in \Kn$, then
$\bx\cdot \by = (x_1y_1\kdots x_ny_n)$. A subgroup $\Ga$ of $\Kn$
is said to be of finite rank if there are $\bu_1\kdots\bu_r\in\Ga$
with the property that for every $\bx\in\Ga$ there are $z\in\Z_{>
0}$ and $z_1,\ldots,z_r\in\Z$ such that $\bx^z
=\bu_1^{z_1}\cdots\bu_r^{z_r}$. The smallest $r$ for which such
$\bu_1\kdots\bu_r$ exist is called the rank of $\Ga$; the rank of
$\Ga$ is equal to $0$ if all elements of $\Ga$ have finite order.

For the moment, let $n=2$. We consider the equation
\begin{equation}
\label{1.3} a_1x_1+a_2x_2=1\quad\mbox{in $\bx =(x_1,x_2)\in\Ga$,}
\end{equation}
where $\ba =(a_1,a_2)\in (K^*)^2$ and where $\Ga$ is a subgroup of
$(K^*)^2$ of finite rank $r$. In 1996, Beukers and Schlickewei
\cite{BS} showed that \formref{1.3} has at most $2^{8(r+2)}$
solutions.

Two pairs $\ba =(a_1,a_2)$, $\bb =(b_1,b_2)$ are called
$\Ga$-equivalent if there is an $\bu\in\Ga$ such that $\bb
=\bu\cdot \ba$. Clearly, two equations \formref{1.3} with
$\Ga$-equivalent pairs of coefficients $\ba$ have the same number
of solutions. In 1988, Gy\H{o}ry, Stewart, Tijdeman and the author
\cite{EGST} showed that there is a finite number of
$\Ga$-equivalence classes, such that for all tuples $\ba
=(a_1,a_2)$ outside the union of these classes, equation
\formref{1.3} has at most \emph{two} solutions.
(In fact they considered only groups $\Ga = U_S\times U_S$
where $U_S$ is the group
of $S$-units in a number field, but their argument works in precisely the
same way for the general case.)
The upper bound
$2$ is best possible. We mention that this result is ineffective
in that the method of proof does not allow to determine the
exceptional equivalence classes. B\'{e}rczes \cite[Lemma 3]{B}
calculated the upper bound $2e^{30^{20}(r+2)}$
for the number of exceptional equivalence classes.

Now let $n\geq 3$. We deal with equations
\begin{equation}
\label{1.1} a_1x_1+\cdots +a_nx_n=1\quad\mbox{in $\bx =(x_1\kdots
x_n)\in\Ga$,}
\end{equation}
where $\ba =(a_1\kdots a_n)\in \Kn$ and where $\Ga$ is a subgroup
of $\Kn$ of finite rank $r$. A solution $\bx$ of \formref{1.1} is
called non-degenerate if
\begin{equation}
\label{1.2} \sum_{i\in I} a_ix_i\not=0\quad\mbox{for each
non-empty subset $I$ of $\{ 1\kdots r\}$.}
\end{equation}
It is easy to show that there are groups $\Ga$ such that any
degenerate solution of \formref{1.1} gives rise to an infinite set
of solutions. Schlickewei, Schmidt and the author \cite{ESS}
showed that equation \formref{1.1} has at most
$e^{(6n)^{3n}(r+1)}$ non-degenerate solutions. Their proof was
based on a version of the quantitative Subspace Theorem, i.e., on
the Thue-Siegel-Roth-Schmidt method. Recently, by a very different
approach based on a method of Vojta and Faltings, R\'{e}mond \cite{Rem}
proved a general quantitative result for subvarieties of tori,
which includes as a special case that for $n\geq 3$
equation \formref{1.1} has at most
$2^{n^{4n^2}(r+1)}$ non-degenerate solutions.

Two tuples $\ba ,\bb \in\Kn$ are called $\Ga$-equivalent if $\bb
=\bu\cdot\ba$ for some $\bu\in\Ga$. Gy\H{o}ry, Stewart, Tijdeman
and the author \cite{EGST} showed that for every sufficiently
large $r$, there are a subgroup $\Ga$ of $(\Q^*)^n$ of rank $r$,
and infinitely many $\Ga$-equivalence classes of tuples $\ba
=(a_1\kdots a_n)\in (\Q^*)^n$, such that equation \formref{1.1}
has at least $e^{2r^{1/2}(\log r)^{-1/2}}$ non-degenerate
solutions. This shows that in contrast to the case $n=2$,
for $n\geq 3$ there is
no uniform bound $C$ independent of $\Ga$ such that for all tuples
$\ba$ outside finitely many $\Ga$-equivalence classes the
number of non-degenerate solutions of \formref{1.1} is at most
$C$.

It turned out to be
more natural to consider the minimal number $m$ such that the set
of solutions of \formref{1.1} can be contained in the union of $m$
proper linear subspaces of $K^n$. Notice that this minimal number
$m$ does not change if $\ba$ is replaced by a $\Ga$-equivalent
tuple. In 1988 Gy\H{o}ry and the author \cite{EG} showed that if $K$ is a
number field and
$\Ga =U_S^n$,
i.e., the $n$-fold direct product of the group of $S$-units in $K$,
then there are finitely
many $\Ga$-equivalence classes $C_1\kdots C_t$ such that for every tuple
$\ba\in\Kn\backslash (C_1\cup\cdots\cup C_t)$ the set of
solutions of \formref{1.1} is contained in the union of not more
than $2^{(n+1)!}$ proper linear subspaces of $K^n$. This was
improved by the author \cite[Thm. 8]{Evb} to $(n!)^{2n+2}$. Both the proofs
of Gy\H{o}ry and the author and that of the author can be extended easily to
arbitrary fields
$K$ of characteristic $0$ and arbitrary subgroups
$\Ga$ of $\Kn$ of finite rank.

For certain special groups $\Ga$, Schlickewei and Viola \cite[Corollary 2]{SV}
improved the author's bound to ${2n+1\choose n}-n^2-n-2$.
In fact, their result is valid for rank one groups
$\Ga = \{ (\alpha_1^z\kdots\alpha_n^z):\, z\in\Z\}$,
where $\alpha_1\kdots\alpha_n$ are non-zero elements of a field $K$
of characteristic $0$ such
that neither $\alpha_1\kdots\alpha_n$, nor any of the quotients
$\alpha_i/\alpha_j$ $(0\leq i<j\leq n)$ is a root of unity.

In the present paper we deduce a further improvement for the general
equation \formref{1.1}.
\\[0.5cm]
{\bf Theorem.} \emph{Let $K$ be a field of characteristic $0$, let
$n\geq 3$, and let $\Ga$ be a subgroup of $\Kn$ of finite rank.
Then there are finitely many
$\Ga$-equivalence classes $C_1\kdots C_t$
of tuples in $\Kn$, such that for every
$\ba =(a_1\kdots a_n)\in \Kn\backslash (C_1\cup\cdots\cup C_t)$,
the
set of non-degenerate solutions of}
\\[0.2cm]
\romref{1.1}\hspace{2.4cm}$a_1x_1+\cdots +a_nx_n=1
\quad\mbox{\emph{in} $\bx =(x_1\kdots x_n)\in\Ga$}$
\\[0.2cm]
\emph{is contained in the union of not more than $2^n$ proper
linear subspaces of $K^n$.}
\\[0.5cm]
We mention that the set of degenerate solutions of \formref{1.1}
is contained in the union of at most $2^n-n-2$ proper linear
subspaces of $K^n$, each defined by a vanishing subsum
$\sum_{i\in I} a_ix_i =0$ where $I$ is a subset of $\{ 1\kdots n\}$
of cardinality $\not= 0,1,n$.
So for
$\ba\not\in C_1\cup\cdots\cup C_t$, the set of (either degenerate
or non-degenerate) solutions of \formref{1.1} is contained in the
union of at most $2^{n+1}-n-2$ proper linear subspaces of
$K^n$.

Our main tool is a qualitative finiteness result due to Laurent
\cite{L} for the number of non-degenerate solutions in $\Ga$ of a
system of polynomial equations  (or rather for the number of
non-degenerate points in $X\cap\Ga$ where $X$ is an algebraic
subvariety of the $n$-dimensional linear torus). Recently,
R\'{e}mond \cite{Rem} established for $K=\overline{\Q}$
an explicit upper bound for the
number of these non-degenerate solutions. Using the latter, it is
possible to compute a (very large) explicit upper bound for the
number $t$ of exceptional equivalence classes, depending on $n$
and the rank $r$ of $\Ga$. We have not worked this out.

In Section 2 we recall Laurent's result. In Section 3 we prove our Theorem.
In Section 4 we give an example showing that our bound $2^n$
cannot be improved to a quantity smaller than $n$.
\\[0.5cm]
\section{Polynomial equations}

Let as before $K$ be a field of characteristic $0$, let $n\geq 2$,
and let $f_1\kdots f_R\in K[X_1\kdots X_n]$ be non-zero polynomials.
Further, let
$\Ga$ be a subgroup of $\Kn$ of finite rank. We consider the
system of equations
\begin{equation}
\label{2.1}
f_i(x_1\kdots x_n )=0\quad (i=1\kdots R)\quad\mbox{in $\bx =(x_1\kdots
x_n)\in\Ga$.}
\end{equation}
Let $\lambda$ be an auxiliary variable.
A solution $\bx =(x_1\kdots x_n)$ of system \formref{2.1} is called \emph{degenerate}
if there are integers $c_1\kdots c_n$ with $\gcd (c_1\kdots c_n)=1$
such that
\begin{equation}
\label{2.2}
f_i(\lambda^{c_1}x_1\kdots \lambda^{c_n}x_n)=0\,\,
\mbox{identically in $\lambda$ for $i=1\kdots R$}
\end{equation}
(meaning that by expanding the expressions, we get linear combinations
of different powers of $\lambda$, all of whose coefficients are $0$).
Otherwise, the solution $\bx$ is called \emph{non-degenerate.}
\\[0.5cm]
{\bf Proposition 2.1.} \emph{System \formref{2.1} has only finitely many
non-degenerate solutions.}
\\[0.2cm]
{\bf Proof.} Without loss of generality we may assume that $K$ is
algebraically closed. Let $X$ denote the set of points $\bx\in
(K^*)^n$ with $f_i(\bx )=0$ for $i=1\kdots R$. By a result of
Laurent \cite[Th\'{e}or\`{e}me 2]{L}, the set of solutions
$\bx\in\Ga$ of \formref{2.1} is contained in the union of finitely
many ``families" $\bx H =\{\bx\cdot \by :\, \by\in H\}$, where
$\bx\in\Ga$ and where $H$ is an irreducible algebraic subgroup of
$\Kn$ such that $\bx H\subset X$. \footnote{For $K=\overline{\Q}$,
R\'{e}mond \cite[Thm. 1]{Rem} showed that the set of solutions of
\formref{2.1} is contained in the union of at most $(nd)^{n^3
m^{3m^2}(r+1)}$ families $\bx H$, where $r$ is the rank of $\Ga$,
$X$ has dimension $m$, and where each polynomial $f_i$ has total
degree $\leq d$. Probably his result can be extended to arbitrary
fields $K$ of characteristic $0$ by means of a specialization
argument.}

Consider a family $\bx H$ with $\bx\in\Ga$, $\bx H\subset X$,
$\dim H >0$. Pick a one-dimensional irreducible algebraic group
$H_0\subset H$. There are integers $c_1\kdots c_n$ with $\gcd
(c_1\kdots c_n)\\=1$ such that $H_0=\{ (\lambda^{c_1}\kdots
\lambda^{c_n}):\,\lambda\in K^*\}$. Then $\bx H_0= \{
(x_0\lambda^{c_0}\kdots x_n\lambda^{c_n}): \,\lambda\in
K^*\}\subset \bx H\subset X$, and the latter implies
\formref{2.2}. Conversely, if $\bx$ satisfies \formref{2.2} then
$\bx H_0\subset X$. Therefore, the solutions of \formref{2.1}
contained in families $\bx H$ with $\dim H >0$ are precisely the
degenerate solutions of \formref{2.1}. Each of the remaining
families $\bx H$, i.e., with $\dim H =0$ consists of a single
solution $\bx$ since $H=\{ (1\kdots 1)\}$. It follows that system
\formref{2.1} has at most finitely many non-degenerate
solutions.\qed
\\[0.5cm]
\section{Proof of the Theorem}

Let again $K$ be a field of characteristic $0$, let $n\geq 3$, and
let $\Ga$ a subgroup of $\Kn$ of finite rank.
Further, let $\ba =(a_1\kdots a_n)\in\Kn$.
We deal with
\\[0.2cm]
\formref{1.1}\hspace{2.4cm}$a_1x_1+\cdots +a_nx_n=1
\quad\mbox{in $\bx =(x_1\kdots x_n)\in\Ga$.}$
\\[0.2cm]
Assume that \formref{1.1} has a non-degenerate solution. By
replacing $\ba$ by a $\Ga$-equivalent tuple we may assume that
$\bone =(1\kdots 1)$ is a non-degenerate solution of
\formref{1.1}. This means that
\begin{equation}
\label{3.1}
\left\{
\begin{array}{l}
a_1+\cdots +a_n=1,\\
\sum_{i\in I} a_i\not= 0\,\,
\mbox{for each non-empty subset $I$ of $\{ 1\kdots  n\}$.}
\end{array}\right.
\end{equation}
We will show that there is a finite set of tuples $\ba$ with
\formref{3.1} such that for each $\ba\in\Kn$ outside this set,
the set of non-degenerate solutions of \formref{1.1} is contained
in the union of not more than $2^n$ proper linear subspaces of $K^n$.
This clearly suffices to prove our Theorem.

By the result of Schlickewei, Schmidt and the author or that
of R\'{e}mond mentioned in
Section 1, there is a finite bound $N$ independent of $\ba$ such
that equation \formref{1.1} has at most $N$ non-degenerate
solutions. (In fact, already Gy\H{o}ry and the author \cite{EG}
proved the existence of such a bound but their method did not
allow to compute it explicitly).

For every tuple $\ba$ with \formref{3.1}, we make a sequence
$\bx_1=\bone$, $\bx_2=(x_{21}\kdots
x_{2n})\kdots$\\$\bx_N=(x_{N1}\kdots x_{Nn})$ such
that each term $\bx_i$ is a non-degenerate solution of
\formref{1.1} and such that each non-degenerate solution of
\formref{1.1} occurs at least once in the sequence.
Then
\\[0.1cm]
\begin{equation}
\label{3.2}
\rank
\left(
\begin{array}{cccc}
1&\cdots&1 &1\\
x_{21}&\cdots &x_{2n}&1\\
\vdots&&\vdots &\vdots\\
\vdots&&\vdots &\vdots\\
x_{N,1}&\cdots &x_{N,n}&1
\end{array}
\right)
\, \leq\, n
\end{equation}
\\[0.1cm]
since
the matrix has $n+1$ linearly dependent columns.
Relation \formref{3.2} means that the determinants of all
$(n+1)\times (n+1)$-submatrices
of the matrix on the left-hand side are $0$.
Thus,
we may view \formref{3.2} as a system of polynomial equations
of the shape \formref{2.1}, to be solved in
$(\bx_2\kdots\bx_N)\in\Ga^{N-1}$. It is important to notice that
this system is independent of $\ba$.

The tuples $\ba$ with \formref{3.1} are now divided into
three classes:

\emph{Class I} consists of those tuples $\ba$ such that
$\rank\{ \bone ,\bx_2\kdots,\bx_N\}=n$ and such that
$(\bx_2\kdots \bx_N)$ is a
non-degenerate solution in $\Ga^{N-1}$ of system \formref{3.2}.

\emph{Class II} consists of those tuples $\ba$ such that
$\rank\{ \bone ,\bx_2\kdots,\bx_N\}<n$.

\emph{Class III} consists of those tuples $\ba$ such that
$(\bx_2\kdots \bx_N)$ is a
degenerate solution in $\Ga^{N-1}$ of system \formref{3.2}.

First let $\ba$ be a tuple of Class I. By
Proposition 2.1, $(\bx_2\kdots \bx_N)$ belongs to a finite set
which is independent of $\ba$. Now $\ba =(a_1\kdots a_n)$ is a
solution of the system of linear equations $a_1+\cdots +a_n=1$,
$x_{i1}a_1+\cdots +x_{in}a_n=1$ $(i=2\kdots N)$. Since by assumption,
$\rank \{ \bone ,\bx_2\kdots \bx_N\}=n$, the tuple
$\ba$ is uniquely determined by $\bx_2\kdots \bx_N$. So
Class I is finite.

For tuples $\ba$ from Class II,
all non-degenerate solutions of \formref{1.1}
lie in a single proper subspace of $K^n$.

Now let $\ba$ be from Class III.
In view of \formref{2.2} this means that
there are integers $c_{ij}$ ($i=2\kdots N$, $j=1\kdots n$), with
$\gcd (c_{ij}:\, i=2\kdots N,\, j=1\kdots n)=1$, such that
\\[0.1cm]
\[
\rank
\left(
\begin{array}{cccc}
1&\cdots&1 &1\\
\lambda^{c_{21}}x_{21}&\cdots &\lambda^{c_{2n}}x_{2n}&1\\
\vdots&&\vdots &\vdots\\
\vdots&&\vdots &\vdots\\
\lambda^{c_{N,1}}x_{N,1}&\cdots &\lambda^{c_{N,n}}x_{N,n}&1
\end{array}
\right)
\, \leq\, n
\]
\\[0.1cm]
identically in $\lambda$, meaning that the determinants of the $(n+1)\times
(n+1)$-submatrices of the left-hand side are identically zero in $\lambda$.

This implies that
there are rational functions $b_j(\lambda )\in K(\lambda )$ $(j=0\kdots n)$,
not all equal to $0$, such that
\begin{equation}
\label{3.3}
\sum_{j=1}^n b_j(\lambda )=b_0(\lambda ) ,\quad
\sum_{j=1}^n b_j(\lambda )\lambda^{c_{ij}}x_{ij}=b_0(\lambda )\,\,\,(i=2\kdots N)\, .
\end{equation}
By clearing denominators, we may assume that
$b_0(\lambda )\kdots b_n(\lambda )$ are polynomials in $K[\lambda ]$
without a common zero.

We substitute $\lambda =-1$. Put $b_j:= b_j(-1)$ $(j=0\kdots n)$ and
$\eps_{ij}:= (-1)^{c_{ij}}$ ($i=2\kdots N$, $j=1\kdots n$).
Then $(b_0\kdots b_n)\not= (0\kdots 0)$, and the numbers $\eps_{ij}$
are not all equal to $1$ since the integers $c_{ij}$ are not all even.
Further, by \formref{3.3} we have
\begin{equation}
\label{3.4}
\left\{
\begin{array}{l}
b_1+\cdots +b_n=b_0\, ,\\
b_1\eps_{i1}x_{i1}+\cdots + b_n\eps_{in}x_{in}=b_0\quad\mbox{for $i=2\kdots N$.}
\end{array}\right.
\end{equation}

We claim that for each tuple $(\eps_1\kdots \eps_n)\in\{ -1,1\}^n$,
the tuple $(b_1\eps_1\kdots b_n\eps_n,b_0)$ is not proportional to
$(a_1\kdots a_n,1)$. Assuming this to be true, it follows from \formref{3.4}
that the set of non-degenerate solutions of \formref{1.1}
is contained in the union of at most $2^n$ proper linear subspaces of $K^n$,
each given by
\[
b_0\big(\sum_{j=1}^n a_jx_j\big)-\sum_{j=1}^n b_j\eps_jx_j=0
\]
for certain $\eps_j\in\{ -1,1\}$ ($j=1\kdots n$).

We prove our claim. First suppose that the tuple $(b_1\kdots b_n,b_0)$
is proportional to
$(a_1\kdots a_n,1)$. There are $i\in\{ 2\kdots N\}$,
$j\in\{ 1\kdots n\}$ such that $\eps_{ij}=-1$. Now $\bx_i$ satisfies both
$\sum_{j=1}^n a_jx_{ij}=1$ (since it is a solution of \formref{1.1})
and $\sum_{j=1}^n a_j\eps_{ij}x_{ij}=1$ (by \formref{3.4}). But then by
subtracting we obtain
$\sum_{j\in J} a_jx_{ij}=0$, where $J$ is the set of indices $j$ with
$\eps_{ij}=-1$. This is impossible since $\bx_i$ is a non-degenerate
solution of \formref{1.1}.

Now suppose that $(b_1\eps_1\kdots b_n\eps_n,b_0)$ is proportional to
$(a_1\kdots a_n,1)$ for certain $\eps_j\in\{ -1,1\}$, not all equal to $1$.
Then by \formref{3.1} and \formref{3.4} we have $\sum_{j=1}^n a_j=1$,
$\sum_{j=1}^n a_j\eps_j=1$.
Again by subtracting, we obtain $\sum_{j\in J} a_j=0$
where $J$ is the set of indices $j$ with $\eps_j=-1$ and this is contradictory
to \formref{3.1}. This proves our claim.

Summarizing, we have proved that Class I is finite, that for every $\ba$
in Class II, all solutions of \formref{1.1} lie in a single proper linear
subspace of $K^n$,
and that for every $\ba$ in Class III, the solutions of \formref{1.1}
lie in the union of $2^n$ proper linear subspaces of $K^n$.
Our Theorem follows.\qed
\\[0.5cm]
\section{Equations whose solutions lie in many subspaces}

We give an example of a group $\Ga$ with the property that there
are infinitely many $\Ga$-equivalence classes of tuples $\ba
=(a_1\kdots a_n)\in\Kn$ such that the set of non-degenerate
solutions of \formref{1.1} cannot be covered by fewer than $n$
proper linear subspaces of $K^n$.

Let $K$ be a field of characteristic $0$, let $n\geq 2$, and let
$\Ga_1$ be an infinite subgroup of $K^*$ of finite rank. Take $\Ga
:= \Ga_1^n=\{ \bx =(x_1\kdots x_n): x_i\in \Ga_1$ for $i=1\kdots
n\}$. Then $\Ga$ is a subgroup of $\Kn$ of finite rank.

Pick $\bu =(u_1\kdots u_n)\in\Ga$ with $b:= u_1+\cdots +u_n\not=0$
and with $\sum_{i\in I} u_i\not=0$ for each non-empty subset $I$ of
$\{ 1\kdots n\}$.
Let $S_n$ denote the group of permutations of $\{ 1\kdots n\}$.
For $\s\in S_n$ write $\bus := (u_{\s (1)}\kdots u_{\s (n)})$.
Then $\bus$ ($\s\in S_n$) are non-degenerate solutions of
\be \label{4.1}
b^{-1}x_1+\cdots + b^{-1}x_n =1\quad\mbox{in $\bx\in\Ga$.}
\ee
For
$i=1\kdots n$, the points $\bus$ with $\s (n) =i$ lie in the
subspace given by
\[
u_i (x_1+\cdots +x_{n-1}) -(b-u_i)x_n=0.
\]
Therefore, for fixed $\bu$, the set $\{ \bus : \s\in S_n\}$ can be
covered by $n$ subspaces. We show that for ``sufficiently general"
$\bu$, this set cannot be covered by fewer than $n$ subspaces.

We need some auxiliary results.
\\[0.5cm]
{\bf Lemma 4.1.} \emph{Let $n\geq 2$ and let $S$ be a subset of
$S_n$ of cardinality $>(n-1)!$. Then there are $\s_1\kdots \s_n\in
S$ such that the polynomial
\be \label{4.2}
F_{\s_1\kdots
\s_n}(X_1\kdots X_n) := \left|
\begin{array}{ccc}
X_{\s_1(1)}&\cdots &X_{\s_1(n)}\\
X_{\s_2(1)}&\cdots &X_{\s_2(n)}\\
\vdots&&\vdots\\
X_{\s_n(1)}&\cdots &X_{\s_n(n)}
\end{array}
\right|
\ee
is not identically zero.}
\\[0.3cm]
{\bf Proof.} We proceed by induction on $n$. For $n=2$ the lemma
is trivial. Assume that $n\geq 3$.

First assume there are $i,j\in\{ 1\kdots n\}$ such that the set
$S_{ij}=\{\s\in S:\,\s (i)=j\}$ has cardinality $>(n-2)!$.
Then after a suitable permutation
of the columns of the determinant of \formref{4.2}
and a permutation of the variables $X_1\kdots X_n$, we obtain that
$S_{nn}$ has cardinality $>(n-2)!$.
The elements of $S_{nn}$
permute $1\kdots n-1$. Therefore, by the induction hypothesis,
there are $\s_1\kdots \s_{n-1}\in S_{nn}$ such that the polynomial
\[
G(X_1\kdots X_{n-1}):=\, \left|
\begin{array}{ccc}
X_{\s_1(1)}&\cdots &X_{\s_1(n-1)}\\
\vdots&&\vdots\\
X_{\s_{n-1}(1)}&\cdots &X_{\s_{n-1}(n-1)}
\end{array}
\right|
\]
is not identically zero. Since $S_{nn}$ has cardinality $\leq
(n-1)!$, there is a $\s_n\in S$ with $\s_n(n)=k\not=n$. Therefore,
\[
F_{\s_1\kdots \s_n}(X_1\kdots X_{n-1},0) =\pm X_k\cdot
G(X_1\kdots X_{n-1})\not =0.
\]
So in particular, $F_{\s_1\kdots \s_n}$ is not identically zero.

Now suppose that for each pair $i,j\in\{ 1\kdots n\}$ the set
$S_{ij}$ has cardinality $\leq (n-2)!$. Together with our
assumption that $S$ has cardinality $>(n-1)!$, this implies that
$S_{ij}\not=\emptyset$ for $i,j\in\{ 1\kdots n\}$. Thus, we may
pick $\s_1\in S$ with $\s_1(1)=1$, $\s_2\in S$ with
$\s_2(2)=1$$\kdots \s_n\in S$ with $\s_n(n)=1$. Then
$F_{\s_1\kdots \s_n}(1,0\kdots 0)=1$, hence
$F_{\s_1\kdots \s_n}$ is not
identically zero.\qed
\\[0.5cm]

Let $T$ denote the collection of tuples $(\s_1\kdots \s_n)$ in
$S_n$ for which $F_{\s_1\kdots \s_n}$ is not identically $0$. Let
$B$ be the set of numbers of the shape $u_1+\cdots +u_n$ where
$\bu =(u_1\kdots u_n)$ runs through all tuples in $\Ga =\Ga_1^n$
with
\be
\label{4.3}
\left\{
\begin{array}{l}
\displaystyle{\sum_{i\in I} u_i\not=0\quad
\mbox{for each $I\subseteq\{1\kdots n\}$ with $I\not=\emptyset$;}}\\
\displaystyle{F_{\s_1\kdots \s_n}(u_1\kdots u_n)\not=0\quad
\mbox{for each $(\s_1\kdots \s_n)\in T$.}}
\end{array}
\right.
\ee
In particular (taking $I=\{ 1\kdots n\}$), each $b\in B$ is non-zero.

Two numbers $b_1,b_2\in K^*$ are called
$\Ga_1$-equivalent if $b_1/b_2\in \Ga_1$.
\\[0.5cm]
{\bf Lemma 4.2.} \emph{The set $B$ is not contained in the union
of finitely many $\Ga_1$-equivalence classes.}
\\[0.3cm]
{\bf Proof.} First suppose that $B\not=\emptyset$.
Assume that $B$ is contained in the union of finitely
many $\Ga_1$-equivalence classes. Let $b_1\kdots b_t$ be
representatives for these classes. Then for every $\bu =(u_1\kdots
u_n)\in\Ga$ with \formref{4.3} there are $b_i\in\{ b_1\kdots
b_t\}$ and $u\in\Ga_1$ such that
\[
u_1+\cdots +u_n =b_iu.
\]
Hence for given $b_i$, $(u_1/u\kdots u_n/u)$ is a non-degenerate
solution of
\[
x_1+\cdots +x_n=b_i\quad\mbox{in $\bx =(x_1\kdots x_n)\in\Ga$.}
\]
Each such equation has only finitely many
non-degenerate solutions. Therefore,
for each $b_i$ there are only
finitely many possibilities for $(u_1/u\kdots u_n/u)$, hence only
finitely many possibilities for $u_1/u_2$.
So if $(u_1\kdots u_n)$ runs through all tuples in $\Ga$ with
\formref{4.3}, then $u_1/u_2$ runs through a finite set, $U$, say.

Now let $F$ be the product of the polynomials
$F_{\s_1\kdots \s_n}$ $((\s_1\kdots \s_n)\in T)$,\\
$\sum_{i\in I} X_i$ $(I\subseteq\{ 1\kdots n\},\, I\not=\emptyset
)$ and $X_1-uX_2$ $(u\in U)$. Then $F(u_1\kdots u_n)=0$ for every
$u_1\kdots u_n\in \Ga_1$. But since $\Ga_1$ is infinite, this
implies that $F$ is identically zero. Thus, if we
assume that $B\not=\emptyset$ and that Lemma
4.2 is false we obtain a contradiction. The assumption $B=\emptyset$ leads
to a contradiction in a similar manner, taking for $F$ the product of the
polynomials $F_{\s_1\kdots \s_n}$ $((\s_1\kdots \s_n)\in T)$,
$\sum_{i\in I} X_i$ $(I\subseteq\{ 1\kdots n\},\, I\not=\emptyset )$.\qed
\\[0.5cm]

Lemma 4.2 implies that the
collection of tuples $(b^{-1}\kdots b^{-1})$ ($n$ times) with
$b\in B$ is not contained in the union of finitely many
$\Ga$-equivalence classes. We show that for every $b\in B$, the
set of non-degenerate solutions of \formref{4.1} cannot be covered
by fewer than $n$ proper linear subspaces of $K^n$.

Choose $b\in B$, and choose $\bu =(u_1\kdots u_n)\in\Ga$ with
$u_1+\cdots +u_n=b$ and with \formref{4.3}. Then each vector
$\bu_{\s}$ $(\s\in S_n)$ is a non-degenerate solution of
\formref{4.1}.

We claim that a proper linear subspace of $K^n$ cannot contain
more than $(n-1)!$ vectors $\bu_{\s}$ ($\s\in S_n$). For suppose
some subspace $L$ of $K^n$ contains more than $(n-1)!$ vectors
$\bu_{\s}$. Then by Lemma 4.1, there are $\s_1\kdots \s_n\in S_n$
such that $\bu_{\s_i}\in L$ for $i=1\kdots n$ and such that
$F_{\s_1\kdots \s_n}$ is not identically $0$.
But since $\bu$ satisfies
\formref{4.3}, we have $F_{\s_1\kdots \s_n}(\bu )\not= 0$.
Therefore, the vectors $\bu_{\s_1}\kdots \bu_{\s_n}$ are linearly
independent. Hence $L=K^n$.

Our claim shows that at least $n$ proper linear subspaces of $K^n$
are needed to cover the set $\bu_{\s}$ $(\s\in S_n)$. Therefore,
the set of non-degenerate solutions of \formref{4.1} cannot lie in
the union of fewer than $n$ proper subspaces.
\\[0.5cm]

\end{document}